\documentclass[12pt]{commalg}

\usepackage{amsmath}
\usepackage{amsthm}
\usepackage[round]{natbib}

\theoremstyle{plain}
\newtheorem{thm}{Theorem}
\newtheorem{lem}[thm]{Lemma}
\newtheorem{prop}[thm]{Proposition}
\newtheorem{cor}[thm]{Corollary}

\newtheorem{ithm}{Theorem}
\newtheorem{icor}[ithm]{Corollary}

\newcommand{\n}{\mathbf{N}}
\newcommand{\z}{\mathbf{Z}}

\DeclareMathOperator{\der}{Der}
\DeclareMathOperator{\depth}{depth}
\DeclareMathOperator{\enm}{End}

\DeclareMathOperator*{\osum}{\oplus}

\begin{document}

\author{Eivind Eriksen}
\title[Connections on modules]{Connections on modules over
quasi-homogeneous plane curves}

\address{Oslo University College \\ Postboks 4 St. Olavs Plass \\
0130 Oslo, Norway} \email{eeriksen@hio.no} \keywords{Connections on
modules, quasi-homogeneous plane curve singularities, simple curve
singularities} \subjclass[2000]{Primary 13N10, 13N15}

\begin{abstract}
Let $k$ be an algebraically closed field of characteristic $0$, and
let $A = k[x,y]/(f)$ be a quasi-homogeneous plane curve. We show
that for any graded torsion free $A$-module $M$ without free
summands, there exists a natural graded integrable connection, i.e.
a graded $A$-linear homomorphism $\nabla: \der_k(A) \to \enm_k(M)$
that satisfies the derivation property and preserves the Lie product.

In particular, a torsion free module $N$ over the complete local
ring $B = \widehat A$ admits a natural integrable connection if $A$
is a simple curve singularity, or if $A$ is irreducible and $N$ is a
gradable module.
\end{abstract}

\maketitle

\section{Introduction}

Let $k$ be an algebraically closed field of characteristic $0$, let
$A$ be a commutative $k$-algebra, and let $M$ be an $A$-module. We
define a \emph{connection} on $M$ to be an $A$-linear homomorphism
    \[ \nabla: \der_k(A) \to \enm_k(M) \]
that satisfy the derivation property, i.e. $\nabla_D(am) = a
\nabla_D(m) + D(a) m$ for all $D \in \der_k(A), \; a \in A$ and $m
\in M$. If $\nabla$ is also a homomorphism of Lie algebras, then we
say that $\nabla$ is an \emph{integrable connection} on $M$.

We say that $A$ is a \emph{quasi-homogeneous plane curve} if $A =
k[x,y]/(f)$ is a reduced, positively graded $k$-algebra given by a
quasi-homogeneous polynomial $f \neq 0$. We shall denote its
completion by $B = \widehat A$, and refer to $B$ as the complete
local ring of a \emph{quasi-homogeneous plane curve singularity}.

The simple curve singularities have been classified by \citet{arn81}
and \citet{wal84}. These singularities have complete local rings of
the form $B = k[[x, y]] / (f)$, where $f$ is one of the following
polynomials:
\begin{align*}
A_n: \; & f = x^2 + y^{n+1} & n \ge 1\\
D_n: \; & f = x^2y + y^{n-1} & n \ge 4\\
E_6: \; & f = x^3 + y^4 & \\
E_7: \; & f = x^3 + xy^3 & \\
E_8: \; & f = x^3 + y^5 &
\end{align*}
In particular, any simple curve singularity is of the form $B =
\widehat A$ for a quasi-homogeneous plane curve $A = k[x,y]/(f)$. By
abuse of language, we shall refer to the graded $k$-algebra $A$ as a
simple curve singularity.

In this paper, we study existence of connections on an $A$-module
$M$ when $A$ is a quasi-homogeneous plane curve and $M$ is a graded
torsion free $A$-module. A connection $\nabla: \der_k(A) \to
\enm_k(M)$ is \emph{graded} if $\nabla_D(M_w) \subseteq
M_{w+\lambda}$ for any homogeneous derivation $D \in \der_k(A)$ of
weight $\lambda$ and any integer $w \in \z$. We show the following
result:

\begin{ithm}
Let $A = k[x,y]/(f)$ be a quasi-homogeneous plane curve, and let $M$
be a graded torsion free $A$-module. If $A$ is irreducible or a
simple curve singularity, then there is a natural graded integrable
connection on $M$.
\end{ithm}

In fact, we believe that any graded torsion free module without free
summands over an arbitrary quasi-homogeneous plane curve admits a
natural graded integrable connection. We give a sufficient condition
in lemma \ref{l:c2-suff}, and show that it is satisfied if $A =
k[x,y]/(y(x^n-y^m))$ for positive integers $m,n \in \n$ with $(m,n)
= 1$ and if $M$ has rank one. However, we have not been able to find
a proof in the general case.

We say that a torsion free $B$-module $N$ is \emph{gradable} if it
is of the form $N \cong \widehat M$ for some graded torsion free
$A$-module $M$. If $B$ is the complete local ring of a simple curve
singularities, then it follows from theorem (15.14) in \citet{yos90}
that any torsion free $B$-module is gradable.

\begin{icor}
Any torsion free module over the complete local ring of a simple
curve singularity admits a natural integrable connection.
\end{icor}

\begin{icor}
Any gradable torsion free module over the complete local ring of an
irreducible quasi-homogeneous plane curve singularity admits a
natural integrable connection.
\end{icor}

\section{Quasi-homogeneous plane curves} \label{s:qhom}

Let $k$ be an algebraically closed field. We define a
\emph{positively graded $k$-algebra} to be a finitely generated
$\n_0$-graded $k$-algebra
    \[ A = \osum_{w \in \n_0} A_w \]
with $A_0 = k$. In this paper, we shall assume that $A$ is a
positively graded, reduced $k$-algebra of the form $A = k[x,y]/(f)$,
where $f \neq 0$. This implies that the complete local ring
$\widehat A$ is a quasi-homogeneous plane curve singularity, and we
shall refer to the graded ring $A$ as a \emph{quasi-homogeneous
plane curve}.

The \emph{simple curve singularities} have been classified by
\citet{arn81} and \citet{wal84}. These singularities have
complete local rings of the form $B = k[[x, y]] / (f)$, where
$f$ is one of the following polynomials:
\begin{align*}
A_n: \; & f = x^2 + y^{n+1} & n \ge 1\\
D_n: \; & f = x^2y + y^{n-1} & n \ge 4\\
E_6: \; & f = x^3 + y^4 & \\
E_7: \; & f = x^3 + xy^3 & \\
E_8: \; & f = x^3 + y^5 &
\end{align*}
In particular, the complete local ring of any simple curve
singularity is of the form $B = k[[x,y]]/(f) \cong \widehat A$,
where $A = k[x,y]/(f)$ is a quasi-homogeneous plane curve. By abuse
of language, we shall refer to the graded $k$-algebra $A$ as a
simple curve singularity.

We introduce some notation that will be used throughout this
paper: Let $k[x,y]$ be a positively graded polynomial
$k$-algebra, with grading given by positive weights $w_x =
\deg(x), \; w_y = \deg(y) \in \n$ satisfying $(w_x, w_y) = 1$,
and let $f$ be a homogeneous element in $k[x,y]$ of positive
weight $w_f = \deg(f) \in \n$. We recall that the homogeneous
elements $f \in k[x,y]$ of weight $w_f$ are the
\emph{quasi-homogeneous polynomials} of weight $w_f$, i.e. the
polynomials $f \in k[x,y]$ that satisfy
    \[ f( \alpha^{w_x} \cdot x, \alpha^{w_y} \cdot y) =
    \alpha^{w_f} \cdot f(x,y ) \]
for all $\alpha \in k^*$. We consider $A = k[x,y]/(f)$ with the
induced grading, which is a positively graded $k$-algebra.

\begin{lem} \label{l:irrhom}
We have that $f = u \cdot f_1 \cdots f_r$, where $u \in k^*$ is
a unit in $k[x,y]$, and $f_i \in k[x,y]$ is an irreducible
polynomial given by $f_i = x$, $f_i = y$ or $f_i = x^{w_y} +
a_i y^{w_x}$ for some $a_i \in k^*$ for $1 \le i \le r$.
\end{lem}
\begin{proof}
In case $w_x = w_y = 1$, it is enough to show that any
quasi-homogeneous polynomial has a factorization into linear
polynomials. This is an elementary fact, since $k$ is
algebraically closed. In the general case, any irreducible
quasi-homogeneous polynomial with more than one summand must
include $x^{s w_y}$ and $y^{s w_x}$. By an argument similar to
the one in case $w_x = w_y = 1$, irreducibility implies that $s
= 1$.
\end{proof}

For $1 \le i \le r$, we denote the weight of $f_i$ by $w_i =
\deg(f_i) \in \n$, and write $A_i = A/(f_i)$ for the
corresponding positively graded domain. We assume that $f_i
\neq f_j$ when $i \neq j$, which means that $A$ is reduced.

For any graded $k$-algebra $R$, let $Q(R)$ be the \emph{graded total
ring of fractions} of $R$, i.e. the localization $Q(R) = S^{-1} R$,
where $S$ is the set a homogeneous non-zero divisors in $R$. Then
the natural map $R \to Q(R)$ is an injective graded
homomorphism, and we identify $R$ with its image in $Q(R)$. We
define the \emph{graded normalization} $\widetilde R$ of $R$ to be
the integral closure of $R$ in $Q(R)$, and write $n: R
\to \widetilde R$ for the graded normalization map.

It follows from lemma \ref{l:irrhom} that there is an
isomorphism $\widetilde A_i \cong k[t_i]$ of graded
$k$-algebras for $1 \le i \le r$, with $d_i = \deg(t_i) \in \n$
given by
    \[ d_i = \begin{cases}
    w_y & \text{if } f_i = x \\
    w_x & \text{if } f_i = y \\
    1 & \text{if } f_i = x^{w_y} + a_i y^{w_x}
    \end{cases} \]
Moreover, the graded normalization map $n_i: A_i \to k[t_i]$ is
the graded algebra homomorphism given by
\begin{align*}
x \mapsto
\begin{cases}
0 & \text{if } f_i = x \\
t_i & \text{if } f_i = y \\
t_i^{w_x} & \text{if } f_i = x^{w_y} + a_i y^{w_x}
\end{cases}
&
&
y \mapsto
\begin{cases}
t_i & \text{if } f_i = x \\
0 & \text{if } f_i = y \\
b_i t_i^{w_y} & \text{if } f_i = x^{w_y} + a_i y^{w_x}
\end{cases}
\end{align*}
for $1 \le i \le r$, where $b_i \in k^*$ is a solution of $a_i
b_i^{w_x} = -1$. Since we have that
    \[ \widetilde A \; \cong \osum_{1 \le i \le r} \widetilde
    A_i, \]
there is an isomorphism $\widetilde A \cong k[t_1] \osum \dots \osum
k[t_r]$ of graded $k$-algebras. We recall that the graded
$k$-algebra $k[t_1] \osum \dots \osum k[t_r]$ has $k$-linear base
$\{  t_i^{l_i} = (0, \dots, t_i^{l_i}, \dots, 0): 1 \le i \le r, \;
l_i \in \n_0 \}$, multiplication given by $t_i^{l_i} \cdot t_j^{m_j}
= \delta_{ij} \cdot t_i^{l_i + m_j}$ for $1 \le i,j \le r, \; l_i,
m_j \in \n_0$, and $k$-algebra structure given by the diagonal map
$k \to k^r$. The grading of $k[t_1] \osum \dots \osum k[t_r]$ is
given by $\deg(t_i^{l_i}) = l_i d_i$ for $1 \le i \le r, \; l_i \in
\n_0$, hence it is a positively graded $k^r$-algebra. We may
identify $A$ with a graded subring of $k[t_1] \osum \dots \osum
k[t_r]$ via the graded normalization map $n: A \to
\widetilde A$.

\section{Derivations}

Let $k$ be an algebraically closed field of characteristic $0$, and
let $R$ be any graded $k$-algebra. We recall that the derivation
module $\der_k(R)$ is the left $R$-module consisting of all
$k$-linear operators $Q: R \to R$ satisfying the Leibniz rule
    \[ Q(r s) = Q(r) s + r \, Q(s) \text{ for all } r,s \in R.
    \]
We say that a derivation $Q \in \der_k(R)$ is homogeneous of weight
$\lambda \in \z$ if $Q(A_w) \subseteq A_{w + \lambda}$ for all $w
\in \n$.

Let $A = k[x,y]/(f)$ be a quasi-homogeneous plane curve. It is
well-known that $\der_k(A)$ is generated, as a left $A$-module,
by the \emph{Euler derivation} $E$ and the \emph{Koszul
derivation} $D$, given by
\begin{align*}
    E & = w_x \, x \cdot \partial_x + w_y \, y \cdot \partial_y
    \\
    D & = f_y \cdot \partial_x - f_x \cdot \partial_y
\end{align*}
We write $\partial_x = \partial / \partial x$ and $\partial_y =
\partial / \partial y$ for the partial derivations, and $f_x =
\partial_x(f)$ and $f_y = \partial_y(f)$ for the partial
derivatives of $f$.

We claim that any derivation $P: A \to A$ has a canonical extension
to a derivation $\widetilde P: \widetilde A \to \widetilde A$. If $r
= 1$, the claim follows from the graded version of a theorem by
\citet{sei66}. If $r > 1$, any derivation $P: A \to A$ induces a
derivation $P_i: A_i \to A_i$ for $1 \le i \le r$, since $E(f_i) =
w_i \cdot f_i$ and $D(f_i) = f_i \cdot [ \partial_x(f_i) \cdot
\partial_y(g_i) - \partial_y(f_i) \cdot \partial_x(g_i) ]$ with
$g_i = f/f_i$, and $\widetilde P = (\widetilde P_1, \dots,
\widetilde P_r): \widetilde A \to \widetilde A$ is an extension of
$P: A \to A$.

\begin{lem}
We have that $\widetilde E = d_1 t_1 \cdot \partial_1 + \dots +
d_r t_r \cdot \partial_r \in \der_k(\widetilde A \,)$, where
$\partial_i = \partial / \partial t_i$ for $1 \le i \le r$.
\end{lem}
\begin{proof}
The Euler derivation $E: A \to A$ is characterized by $E(h) = w
\cdot h$ for any homogeneous element $h \in A$ of weight $w$.
Since $n: A \to \widetilde A$ is a graded algebra homomorphism,
it follows that $\widetilde E: \widetilde A \to \widetilde A$
can be characterized in the same way, and $\widetilde E = d_1
t_1 \cdot
\partial_1 + \dots + d_r t_r \cdot \partial_r$.
\end{proof}

Let $\Gamma_i \subseteq \n_0$ be the semigroup given by
$\Gamma_i = \{ \gamma \in \n_0: t_i^{\gamma} \in A \}$ for $1
\le i \le r$, where we identify $A$ with its image in
$\widetilde A$. Hence we have
    \[ \gamma \in \Gamma_i \Longleftrightarrow (0, \dots,
    t_i^{\gamma}, \dots, 0) \in A \]
for any $\gamma \in \n_0$. We remark that if $r > 1$, then $0
\not \in \Gamma_i$ for $1 \le i \le r$, so $\Gamma_i$ should be
thought of as a semigroup without zero in that case. It is
clear that $\n_0 \setminus \Gamma_i$ is finite for $1 \le i \le
r$, hence we may define $c_i = \min \{ \gamma \in \Gamma_i:
\gamma + \n_0 \subseteq \Gamma_i \}$ to be the conductor of
$\Gamma_i$, and $g_i = c_i -1$ to be the Frobenius number of
$\Gamma_i$.

\begin{prop}
We have that $\widetilde D = \beta_1 t_1^{c_1} \cdot \partial_1
+ \dots + \beta_r t_r^{c_r} \cdot \partial_r  \in
\der_k(\widetilde A \,)$ for some constants $\beta_1, \dots,
\beta_r \in k^*$, where $\partial_i = \partial / \partial t_i$
for $1 \le i \le r$. Moreover, the conductor $c_i$ of
$\Gamma_i$ is given by
    \[ c_i = \frac{w_f-w_i}{d_i} + c(A_i) \]
for $1 \le i \le r$, where $c(A_i)$ is the conductor of the
irreducible quasi-homogeneous curve $A_i$.
\end{prop}
\begin{proof}
If $r = 1$, then it is well-known that $\widetilde D = \beta_1
t_1^{c_1} \cdot \partial_1$ for some $\beta_1 \in k^*$, see for
instance section 6 of \citet{er03}. Hence we may assume that $r
> 1$. Let $g_i = f/f_i \in k[x,y]$ for $1 \le i \le r$. It is
clear that $n_I(g_i) = 0$ if $I \neq i$ and that $n_i(g_i) \neq 0$.
In fact, we see that $n_i(g_i) \in k[t_i]$ is homogeneous of degree
$w - w_i$, and therefore of the form
    \[ n_i(g_i) = \beta_i' \cdot t_i^{\frac{w_f-w_i}{d_i}} \]
for some $\beta_i' \in k^*$, hence $\frac{w - w_i}{d_i} \in
\Gamma_i$. For $1 \le i \le r$, we also see that any element $h_i
\in k[x,y]$ with the property that $n_I(h_i) = 0$ for $I \neq i$
must be of the form $h_i = q \cdot g_i$ with $q \in k[x,y]$.
Therefore, $c_i = \frac{w - w_i}{d_i} + c(A_i)$, where $c(A_i)$ is
the conductor of $A_i$. Explicitly, $c(A_i) = 0$ if $f_i = x$ or
$f_i = y$, and $c(A_i) = (w_x -1)(w_y - 1)$ otherwise.

On the other hand, we see that the derivation $D_i \in
\der_k(A_i)$ induced by the Koszul derivation $D \in \der_k(A)$
has the form
    \[ D_i = g_i \cdot ( \partial_y(f_i) \, \partial_x -
    \partial_x(f_i) \, \partial_y ) \]
for $1 \le i \le r$. Since $\partial_y(f_i) \, \partial_x -
\partial_x(f_i) \, \partial_y$ is the Koszul derivation on
$A_i$, its natural extension to $\widetilde A_i$ has the form
$\beta_i'' \, t_i^{c(A_i)} \cdot \partial_i$ for some constant
$\beta_i'' \in k^*$, see for instance section 6 of
\citet{er03}. It follows that $\widetilde D_i \in
\der_k(\widetilde A_i)$ has the form
    \[ \widetilde D_i = n_i(g_i) \cdot \beta_i'' \, t_i^{c(A_i)}
    \cdot \partial_i = \beta_i \, t^{c_i} \cdot \partial_i \]
with $\beta_i = \beta_i' \cdot \beta_i'' \in k^*$.
\end{proof}

We remark that $\widetilde D =  q \cdot  \widetilde E$ considered as
elements in $\der_k(\widetilde A)$, where $q \in \widetilde A$ is
given by $q = \bigl( \frac{\beta_1}{d_1} t_1^{g_1}, \dots,
\frac{\beta_r}{d_r} t_r^{g_r} \bigr)$. In fact, $q \in (A:m)$, where
$m = (x,y) \subseteq A$ denotes the maximal homogeneous ideal of
$A$, and one may show that $\der_k(A) = (A:m) \cdot E$.

\section{Graded torsion free modules} \label{s:tfmod}

Let $k$ be an algebraically closed field, let $A = k[x,y]/(f)$
be a quasi-homogeneous plane curve, and let $M$ be a finitely
generated, graded $A$-module. If the natural map $M \to M
\otimes_A Q(A)$ is injective, we say that $M$ is \emph{torsion
free}. Notice that $M$ is torsion free if and only if it is
maximal Cohen-Macaulay, i.e. $\depth M_m = 1$, where $m = (x,y)
\subseteq A$ is the maximal homogeneous ideal of $A$.

\begin{lem}
If $M$ is a torsion free $A$-module, then we may identify $M$
with a graded submodule of
    \[ \osum_{1 \le i \le r} \widetilde A_i \, ^{s_i} \cong k[t_1]^{s_1}
    \oplus \dots \oplus k[t_r]^{s_r}, \]
where $s_i \ge 0$ is the rank of $M \otimes_A A_i$ as an
$A_i$-module for $1 \le i \le r$.
\end{lem}
\begin{proof}
It is clear that $M \otimes_A A_i$ has a rank $s_i \ge 0$ as an
$A_i$-module for $1 \le i \le r$, since $Q(A_i) \cong
k[t_i,t_i^{-1}]$. We have
    \[ M \otimes_A Q(A) \cong \osum_{1 \le i \le r} M \otimes_A
    Q(A_i), \]
and $M \otimes_A Q(A_i) \cong (M \otimes_A A_i) \otimes_{A_i} Q(A_i)
\cong Q(A_i)^{s_i}$ for $1 \le i \le r$. Hence the natural map $M
\hookrightarrow M \otimes_A Q(A) \cong \osum_i Q(A_i)^{s_i}$
identifies $M$ with a graded submodule of $\osum_i Q(A_i)^{s_i}$.
Since $M$ is finitely generated, we may identify $M$ with a graded
submodule of $\osum_i \widetilde A_i \, ^{s_i} = \osum_i
k[t_i]^{s_i}$ by multiplying with a suitable denominator.
\end{proof}

For $1 \le i \le r$, we choose a set $\{ e_{ij}: 1 \le j \le
s_i \}$ of homogeneous generators for the graded $\widetilde
A_i$-module $\widetilde A_i \, ^{s_i}$, and write $f_{ij} =
\deg(e_{ij}) \in \z$ for their weights. We consider $\{ e_{ij}:
1 \le i \le r, \; 1 \le j \le s_i \}$ as a set of homogeneous
generators for the graded $\widetilde A$-module $\osum_i
\widetilde A_i \, ^{s_i}$, where $\widetilde A$ acts on
$\osum_i \widetilde A_i \, ^{s_i}$ via the projections $\pi_i:
\widetilde A \to \widetilde A_i$ for $1 \le i \le r$.

Let $M \subseteq \osum_i \widetilde A_i \, ^{s_i}$ be a graded
submodule, let $\{ m_l: 1 \le l \le L \}$ be a set of homogeneous
generators of $M$, and write $w_l = \deg(m_l) \in \z$ for their
weights. Explicitly, $m_l$ is of the form
    \[ m_l = \sum_{i,j} a_{ij}(l) \cdot e_{ij}, \]
for $1 \le l \le L$, where $a_{ij}(l) \in \widetilde A_i$ is
homogeneous of weight $w_l - f_{ij}$ for $1 \le i \le r$, $1 \le j
\le s_i$. Let us write $\pi_i: \osum_i \widetilde A_i \, ^{s_i} \to
\widetilde A_i \; ^{s_i}$ for the $i$'th projection. We may assume
that the generating set $\{ m_l: 1 \le l \le L \}$ satisfy the
following condition:
\begin{enumerate} \label{c:one}
    \item[(C1)] For all $1 \le i \le r, \; 1 \le j \le
        s_i$, there is a homogeneous generator $m_l$ with
        $\pi_i(m_l) = e_{ij}$
\end{enumerate}
In particular, condition (C1) implies that we may choose set of
homogeneous generators containing $\{ e_1, \dots, e_s \}$ when $r =
1$.

\section{Connections on graded modules}

Let $k$ be an algebraically closed field of characteristic $0$, and
let $R$ be any graded $k$-algebra. We define a \emph{connection} on
a graded $R$-module $V$ to be an $R$-linear homomorphism
    \[ \nabla: \der_k(R) \to \enm_k(V) \]
that satisfies the derivation property, i.e. $\nabla_D(rv) = r
\nabla_D(v) + D(r) v$ for all $D \in \der_k(R), \; r \in R, \; v \in
V$. We say that $\nabla$ is \emph{graded}, or homogeneous of weight
$0$, if $\nabla_D(V_w) \subseteq V_{w + \lambda}$ for any
homogeneous derivation $D \in \der_k(R)$ of weight $\lambda$ and for
any integer $w \in \z$, and that $\nabla$ is \emph{integrable} if it
is a Lie algebra homomorphism.

Let $A = k[x,y]/(f)$ be a quasi-homogeneous plane curve, and let $M$
be a finitely generated, graded torsion free $A$-module. We may
assume that $M \subseteq \osum_i \widetilde A_i \, ^{s_i}$ is a
graded submodule that satisfies (C1).

The Euler derivation $E: A \to A$ has a canonical extension to a
derivation $\widetilde E = (\widetilde E_1, \dots, \widetilde E_r):
\widetilde A \to \widetilde A$, and an induced $k$-linear action
$\nabla_E \in \enm_k(\osum_i \widetilde A_i \, ^{s_i})$ that
satisfies the derivation property. The induced action $\nabla_E$ is
characterized by $\nabla_E(h) = w \cdot h$ for any homogeneous
element $h \in \osum_i \widetilde A_i \, ^{s_i}$ of weight $w$.

The Koszul derivation $D: A \to A$ has a canonical extension to a
derivation $\widetilde D = (\widetilde D_1, \dots, \widetilde D_r):
\widetilde A \to \widetilde A$, and $\widetilde D = q \cdot
\widetilde E$ in $\der_k(\widetilde A)$ with $q \in (A:m)$. Hence
there is a unique extension of the $k$-linear action $\nabla_E$ to a
connection $\nabla: \der_k(A) \to \enm_k(M)$, with induced
$k$-linear action $\nabla_D \in \enm_k(\osum_i \widetilde A_i \,
^{s_i})$ given by $\nabla_D(h) = q \cdot \nabla_E(h) = q w \cdot h$
for any homogeneous element $h \in \osum_i \widetilde A_i \, ^{s_i}$
of weight $w$. We remark that $\nabla$ is a graded integrable
connection.

\begin{lem} \label{l:koszul}
Let $M \subseteq \osum_i \widetilde A_i \, ^{s_i}$ be a graded
submodule. If $\nabla_D(M) \subseteq M$, then $\nabla$ induces a
graded integrable connection on $M$.
\end{lem}

For a graded submodule $M \subseteq \osum_i \widetilde A_i \,
^{s_i}$ that satisfies (C1), we consider the following conditions:
\begin{enumerate}
    \item[(C2)] For $1 \le i \le r, \; 1 \le j \le s_i$, we have
    $t_i^{g_i} \, e_{ij} \in M$
    \item[(C3)] For $1 \le i \le r, \; 1 \le j \le s_i$, we have
    $f_{ij} = \lambda$ for some $\lambda \in \z$.
\end{enumerate}
We remark that given a description of $M \subseteq \osum_i
\widetilde A_i \, ^{s_i}$ in concrete terms, it is easy to check if
the conditions (C2) and (C3) hold.

\begin{lem} \label{l:c2-suff}
Let $M \subseteq \osum_i \widetilde A_i \, ^{s_i}$ be a graded
submodule. If $M$ satisfies (C1) and (C2), then $\nabla_D(M)
\subseteq M$.
\end{lem}
\begin{proof}
If $\gamma \in \Gamma_i$, then $t_i^{\gamma} \in A$ by definition,
hence condition (C1) implies that $t_i^{\gamma} e_{ij} \in M$. In
particular, $t_i^{\gamma} e_{ij} \in M$ for all $\gamma \ge g_i$ by
(C1) and (C2). On the other hand, the homogeneous generator $m_l$ of
$M$ of weight $w_l$ has the form $m_l = \sum_{ij} a_{ij}(l) \cdot
e_{ij}$ for $1 \le l \le L$, where $a_{ij}(l) \in \widetilde A_i$ is
homogeneous of weight $w_l - f_{ij} \ge 0$ for $1 \le i \le r, \; 1
\le j \le s_i$. Hence we have that $\nabla_D(m_l) = w_l \, q \cdot
m_l \in M$ for $1 \le l \le L$, and it follows that $\nabla_D(M)
\subseteq M$.
\end{proof}

\begin{lem} \label{l:equal-suff}
Let $M \subseteq \osum_i \widetilde A_i \, ^{s_i}$ be a graded
submodule. If $M$ satisfies (C1) and (C3), then
$\nabla_D(M(-\lambda)) \subseteq M(-\lambda)$.
\end{lem}
\begin{proof}
We may assume that $\lambda = 0$, hence $w_l = \deg(m_l) \ge 0$ for
all $l$. If $w_l = 0$, then $\deg(a_{ij}(l)) = w_l - f_{ij} = 0$ for
all $i,j$, hence $\nabla_D(m_l) = 0$. If $w_l > 0$, then
$\deg(a_{ij}(l)) = w_l - f_{ij} = w_l > 0$ for all $i,j$, hence $q
\, a_{ij}(l) \in A$ and $\nabla_D(m_l) = \sum_{ij} w_l \, q \,
a_{ij}(l) \cdot e_{ij} \in M$.
\end{proof}

In particular, if $M$ is free and indecomposable, then we may
identify $M$ with the submodule of $\widetilde A = \osum_i
\widetilde A_i$ generated by $e_{11} + e_{21} + \dots + e_{r1}$.
Hence (C3) holds with $\lambda = f_{11} = \dots = f_{r1}$, and there
is a natural integrable graded connection on $M(-\lambda)$ by lemma
\ref{l:koszul} and lemma \ref{l:equal-suff}.

\begin{prop} \label{p:irr}
If $A$ is irreducible and $M$ is indecomposable and not free, then
$M$ satisfies (C1) and (C2).
\end{prop}
\begin{proof}
Let $M \subseteq \widetilde A \, ^s$ be a graded submodule that
satisfies (C1). We may assume that the set $\{ m_l: 1 \le l \le L
\}$ of homogeneous generators of $M$ satisfies $L > s$ and $m_l =
e_l$ for $1 \le l \le s$. For $l > s$, $m_l \in M$ has the form
    \[ m_l = \sum_{1 \le j \le s} a_j(l) \cdot e_j \]
where $a_j(l) \in \widetilde A$ is homogeneous of weight
$\deg(a_j(l)) = w_l - f_j$ for $1 \le j \le s$. We may assume that
$a_j(l) = 0$ when $w_l - f_j \in \Gamma$, since $e_j \in M$. We may
also assume that $a_j(l) = 0$ or $a_J(l) = 0$ when $j \neq J$ and
$\deg(a_J(l)) - \deg(a_j(l)) \in \Gamma$, since $M$ is
indecomposable.

We claim that these assumptions imply that $t^g e_j \in M$ for $1
\le j \le s$. In fact, $a_j(l) \neq 0$ for some integer $l>s$,
otherwise $A \cdot e_j$ would be a direct summand of $M$. By
assumption, $\deg(a_j(l)) \not \in \Gamma$, and since $\Gamma$ is
symmetric, it follows that $g - \deg(a_j(l)) \in \Gamma$ and
therefore that $b = t^{g - \deg(a_j(l))} \in A$. Hence $b \cdot m_l
\in M$ can be expressed as
    \[ b \cdot m_l = c \, t^g \cdot e_j + \sum_{J \neq j} b \,
    a_J(l) \cdot e_{J} \]
for some $c \in k^*$. If $a_J(l) \neq 0$ for some $J \neq j$, then
by assumption, $\deg(a_j(l)) - \deg(a_J(l)) \not \in \Gamma$, and $g
- ( \deg(a_j(l)) - \deg(a_J(l)) ) \in \Gamma$. It follows that
$\deg(b \, a_J(l)) \in \Gamma$ and that $b \, a_J(l) \in A$ when $J
\neq j$ and $a_J(l) \neq 0$. Hence $b \, a_J(l) \cdot e_J \in M$ for
all $J \neq j$, and this implies that $t^g e_j \in M$.
\end{proof}

It seems probable that any indecomposable module $M$ that is not
free satisfies (C1) and (C2), even when $A$ has more than one
irreducible component, but we have not been able to find a proof in
the general case. To illustrate the combinatorial problems involved,
we give an example. This example also shows that there are modules of
rank one that do not satisfy (C3).

Let $A = k[x,y]/(y(x^n-y^m))$ for some natural numbers $m,n \in \n$
with $(m,n) = 1$. Then $A \cong k[(t_1,t_2^m), (0,t_2^n)]\subseteq
k[t_1] \osum k[t_2]$ and the semigroups $\Gamma_i$ are given by
\begin{align*}
    \Gamma_1 & = \{ n, n+1, n+2, \dots \} \\
    \Gamma_2 & = \langle n, n+m, n+2m, \dots, n + (n-1)m \rangle
\end{align*}
with $g_1 = n-1$ and $g_2 = mn-m = m(n-1)$. Let $M$ be a finitely
generated graded torsion free $A$-module of rank one, and assume
that $M$ is not free. Then we may assume that $M$ has a set of
homogeneous generators containing one of the following sets:
\begin{enumerate}
    \item $\{ e_{11} + e_{21}, t_2^h \cdot e_{21} \}$ for $h \in \n
    \setminus \Gamma_2$
    \item $\{ e_{11} + t_2^h \cdot e_{21}, e_{21} \}$ for $h \in \n
    \setminus \langle m,n \rangle$
\end{enumerate}
In the first case, we see that $g_2 - h \in \langle m,n \rangle$,
hence there is an element $a = (*,t_2^{g_2 - h}) \in A$ and $a \cdot
t_2^h \, e_{21} = t_2^{g_2} \, e_{21} \in M$. In fact, $h - n \not
\in \langle m,n \rangle$ since $h \not \in \Gamma_2$, hence $g(A_2)
- (h - n) = g_2 - h \in \langle m,n \rangle$. Moreover,
$(t_1,t_2^m)^{n-1} \cdot (e_{11} + e_{21}) = t_1^{g_1} e_{11} +
t_2^{g_2} e_{21} \in M$. It follows that condition (C2) holds. In
the second case, $(t_1,t_2^m)^{n-1} \cdot e_{21} = t_2^{g_2} e_{21}
\in M$, and $(t_1,t_2^m)^{n-1} \cdot (e_{11} + t^h e_{21}) =
t_1^{g_1} e_{11} + t_2^{g_2 + h} e_{21} \in M$. Since $g_2 + h \in
\Gamma_2$, it follows that condition (C2) holds. We conclude that
(C2) holds in both cases, while (C3) holds only in the first case.

\begin{thm}
Let $A = k[x,y]/(f)$ be a quasi-homogeneous plane curve, and let $M$
be a finitely generated graded torsion free $A$-module. If $A$ is
irreducible or a simple curve singularity, then there is a natural
graded integrable connection on $M$.
\end{thm}
\begin{proof}
If $M$ is a free and indecomposable, then $M$ admits a natural
integrable connection by the comment following lemma
\ref{l:equal-suff}. By lemma \ref{l:koszul} and \ref{l:c2-suff}, it
is therefore enough to show that there is a graded embedding $M
\subseteq \osum_i \widetilde A_i \, ^{s_i}$ such that conditions
(C1) - (C2) hold when $M$ is indecomposable and non-free. Hence the
theorem follows from proposition \ref{p:irr} when $A$ is
irreducible. When $A$ is a simple curve singularity with more than
one irreducible component, it is straight-forward to verify that
(C2) holds in each case, using the complete list of indecomposable
modules given in the appendix of \citet{gr-kn85}.
\end{proof}

\section{Connections on modules over curve singularities}

Let $k$ be an algebraically closed field of characteristic $0$, and
let $A = k[x,y]/(f)$ be a quasi-homogeneous plane curve. We denote
its $m$-adic completion by $B = \widehat A \cong k[[x,y]]/(f)$,
where $m = (x,y) \subseteq A$ is the maximal homogeneous ideal of
$A$. Hence $B$ is the complete local ring of a quasi-homogeneous
plane curve singularity.

If $M$ is a graded torsion free $A$-module, then $\widehat M$ is a
torsion free $B$-module. We say that a torsion free $B$-module $N$
is \emph{gradable} if there is a graded torsion free $A$-module $M$
such that $\widehat M \cong N$. We remark that by theorem (15.14) in
\citet{yos90} and the fact that the simple curve singularities have
finite CM representation type, it follows that any torsion free
module over a simple curve singularity is gradable.

\begin{cor} \label{c:simple}
Any torsion free module over the complete local ring of a simple curve
singularity admits a natural integrable connection.
\end{cor}

\begin{cor}
Any gradable torsion free module over the complete local ring of an
irreducible quasi-homogeneous plane curve singularity admits a
natural integrable connection.
\end{cor}

\section*{acknowledgements}

The results concerning simple curve singularities were obtained
while I worked on my Master's thesis \cite{er94} at the University
of Oslo. I thank my thesis supervisor, professor O.A.
Laudal, for introducing me to this problem and for sharing his
ideas. I also thank the referee for his careful reading of the 
manuscript and his suggestions for improving the exposition.

\bibliographystyle{commalg}
\bibliography{eeriksen}

\begin{thebibliography}{7}
\providecommand{\natexlab}[1]{#1}
\expandafter\ifx\csname urlstyle\endcsname\relax
  \providecommand{\doi}[1]{doi:\discretionary{}{}{}#1}\else
  \providecommand{\doi}{doi:\discretionary{}{}{}\begingroup
  \urlstyle{rm}\Url}\fi

\bibitem[{Arnol{\cprime}d(1981)}]{arn81}
Arnol{\cprime}d, V.~I. (1981).
\newblock \emph{Singularity theory}, volume~53 of \emph{London Math. Soc.
  Lecture Note Ser.}
\newblock Cambridge Univ. Press.

\bibitem[{Eriksen(1994)}]{er94}
Eriksen, E. (1994).
\newblock \emph{Konneksjoner p{\aa} en klasse moduler over simple
  kurve\-singulariteter}.
\newblock Master's thesis, University of Oslo.

\bibitem[{Eriksen(2003)}]{er03}
Eriksen, E. (2003).
\newblock Differential operators on monomial curves.
\newblock \emph{J. Algebra}, 264(1):186--198.

\bibitem[{Greuel and Kn\"orrer(1985)}]{gr-kn85}
Greuel, G.-M. and Kn\"orrer, H. (1985).
\newblock Einfache {K}urvensingularit\"aten und torsionsfreie {M}oduln.
\newblock \emph{Math. Ann.}, 270:417--425.

\bibitem[{Seidenberg(1966)}]{sei66}
Seidenberg, A. (1966).
\newblock Derivations and integral closure.
\newblock \emph{Pacific J. Math.}, 16:167--173.

\bibitem[{Wall(1984)}]{wal84}
Wall, C. T.~C. (1984).
\newblock Notes on the classification of singularities.
\newblock \emph{Proc. London Math. Soc. (3)}, 48(3):461--513.

\bibitem[{Yoshino(1990)}]{yos90}
Yoshino, Y. (1990).
\newblock \emph{Cohen-{M}acaulay modules over {C}ohen-{M}acaulay rings}, volume
  146 of \emph{London Math. Soc. Lecture Note Ser.}
\newblock Cambridge Univ. Press.

\end{thebibliography}

\end{document}